\documentclass[11pt]{revtex4-1}
\usepackage[utf8]{inputenc}
\usepackage[T1]{fontenc}
\usepackage{graphicx}
\usepackage{soul}
\usepackage{amssymb}
\usepackage{amsmath}
\usepackage{amsthm}
\usepackage{mathtools}
\usepackage[dvipsnames]{xcolor}
\usepackage[colorlinks=false,pdfborder={0 0 0}]{hyperref}
\usepackage{microtype}
\usepackage{cleveref}
\usepackage{url}
\usepackage{dcolumn}
\usepackage{pdfpages}
\usepackage[overload]{textcase}

\colorlet{RED}{red}
\title{}
\date{}
\hypersetup{
  pdfkeywords={},
  pdfsubject={},
  pdfcreator={Emacs Org-mode version 7.9.3f}}

\begin{document}

\newtheorem{theorem}{Theorem}
\title{Non-self-similar blow-up in the heat flow for harmonic maps in higher dimensions}
\begin{abstract}
  We analyze the finite-time blow-up of solutions of the heat flow for
  $k$-corotational maps $\mathbb R^d\to S^d$.  For each dimension
  $d>2+k(2+2\sqrt{2})$ we construct a countable family of blow-up
  solutions via a method of matched asymptotics by glueing a re-scaled
  harmonic map to the singular self-similar solution: the equatorial
  map.  We find that the blow-up rates of the constructed solutions
  are closely related to the eigenvalues of the self-similar
  solution.  In the case of $1$-corotational maps our solutions are
  stable and represent the generic blow-up.
\end{abstract}

\author{Paweł Biernat}
\email{pawel.biernat@gmail.com}
\affiliation{Institute of Mathematics, Jagiellonian University, Kraków}
\date{\today}

\keywords{type II, matched asymptotics, singularity, parabolic,
harmonic map flow, asymptotics}

\maketitle

\section*{Introduction}
\label{sec-1}
\subsection*{}

\label{sec-1-1}

    A map $F:M\to N\subset \mathbb R^k$ between two compact Riemannian
    manifolds $M$ and $N$ is called harmonic if it is a critical point
    of the functional
    \begin{align*}
    E(F)=\frac{1}{2}\int_M |\nabla F|^2\,dV_M.
    \end{align*}
    The heat flow for harmonic maps was introduced by Eells and Sampson
    \cite{Eells1964} as a method of deforming any smooth map $F_0$ to a
    harmonic map via the equation
    \begin{align}
    \label{eq:80}
    \partial_t F=(\Delta F)^\top,\quad F|_{t=0}=F_0,
    \end{align}
    where $(\Delta F)^\top$ is a projection of $(\Delta F)\in \mathbb
    R^k$ to $T_{F}N$---a tangent space to $N$ at the point $F$.  For
    any solution to \eqref{eq:80} we have
    \begin{align*}
    \frac{d}{dt}E(F)=-\int_M|\partial_t F|^2\, dV_M\le 0.
    \end{align*}
    If the flow exists for all times, $E(F)\ge0$ converges to some
    $E_\infty$, suggesting that $F\to F_\infty$ with $F_\infty$ being
    a harmonic map.  This approach proved to work only for target
    manifolds $N$ with non-positive sectional curvature.  If there is
    a point in $N$ with positive sectional curvature, the gradient of
    a solution to \eqref{eq:80} may blow-up in a finite time.  In
    consequence, existence of global in time solutions may be
    established only in a weak sense \cite{Chen1989}.  Moreover, the
    uniqueness of solutions can no longer be guaranteed
    \cite{Chen1989}.  For explicit examples of non-unique weak
    solutions to \eqref{eq:80} in the case of maps $\mathbb R^d\to
    S^d$ with $3\le d\le6$ see \cite{Biernat2011} and
    \cite{Germain2010}.

\label{sec-1-2}

    In order to overcome the problems posed by a finite-time blow-up
    and to investigate the circumstances in which the uniqueness is
    lost one has to fully understand the blow-up mechanism.  The most
    general classification of solutions with a blow-up divides them
    into two types.  We call a solution $F$ to \eqref{eq:80} that
    blows up in finite time $T$ to be of type I if there exists a
    constant $C$ such that
    \begin{align}
    \label{eq:81}
    (T-t) \sup_M|\nabla F|^2 \le C
    \end{align}
    holds for $t<T$ where $T$ is the blow-up time; if \eqref{eq:81}
    does not hold the blow-up is of type II.

\label{sec-1-3}

    The reason for this classification becomes clear when we take maps
    $\mathbb R^d\to N$.  Then, if the blow-up is of type I, we know
    that $F(x,t)=w\left(\frac{x-x_0}{\sqrt{T-t}}\right)$ near an
    isolated singularity located at $(x_0,T)$
    \cite[p.~293]{Struwe1996}.  The function $w:\mathbb R^d\to N$
    describes the profile of a singular solution $F$ and the question
    of existence of singular solutions of type I reduces to the
    question of existence of admissible profile functions $w$.  When
    the blow-up is of type II there is no similar universal
    description of what $F$ looks like near the singularity and any
    type II solution has to be considered on a case by case basis.

    A careful reader will notice that $\mathbb R^d$ is not a compact
    manifold.  Because in this paper we consider only isolated
    singularities it is a matter of convenience to replace the compact
    domain $M$ with a non-compact tangent space $T_{x_0}M=\mathbb
    R^d$, i.e. to neglect the curvature of the domain.  Such
    simplification does not affect the blow-up mechanism.

\label{sec-1-4}

    Let us consider the simplest positively curved target, $S^d$
    embedded in $\mathbb R^{d+1}$ in a canonical way.  The deformation
    of a map $\mathbb R^{d}\to S^{d}$ according to the harmonic map
    heat flow \eqref{eq:80} simplifies to
    \begin{align}
    \label{eq:83}
    \partial_t F=\Delta F+|\nabla F|^2 F.
    \end{align}
    Let us introduce spherical coordinates $(r,\omega)$ on $\mathbb
    R^d$ and coordinates $(u,\Omega)$ on $S^d$, with $u$ denoting the
    latitudinal position on $S^d$ and $\Omega\in S^{d-1}$
    parametrizing the equator.  Using these coordinates we can further
    restrict $F$ to a highly symmetric class of $k$-corotational maps
    \begin{equation}
    \label{eq:F-ansatz}
    (r,\omega)\to(u(r,t),\Omega_k(\omega)).
    \end{equation}
    $\Omega_k$ is a (non-constant) harmonic map with a constant energy
    density $\lvert \nabla \Omega_k\rvert^2=k(k+d-2)$, the number
    $k=1,2,3,\dots$ corresponds to a topological degree of map
    \eqref{eq:F-ansatz}.  The class of $k$-corotational maps is
    preserved by the harmonic map flow and the ansatz
    \eqref{eq:F-ansatz} reduces \eqref{eq:83} to
    \begin{align}
    \label{eq:u}
    \partial_t u=\frac{1}{r^{d-1}}\partial_r\left(r^{d-1}\partial_r
    u\right)-\frac{k(d+k-2)}{2r^2}\sin(2u).
    \end{align}
    The Dirichlet energy $E(F)$ can be expressed (up to a
    multiplicative constant) in terms of $u$ as
    \begin{equation}
    \label{eq:Eu}
    E(u)=\frac{1}{2}\int_0^\infty \left((\partial_r u)^2+k(d+k-2)\frac{\sin^2(u)}{r^2}\right)r^{d-1}\,dr
    \end{equation}
    Regularity of $F$ enforces a boundary condition $u(0,t)=0$, while
    boundary condition at $r=\infty$ follows from the finiteness of
    Dirichlet energy $E(u)<\infty$.  The monotonicity of energy
    \begin{equation}
    \label{eq:Eudt}
    \frac{d}{dt}E(u)=-\int_0^\infty (\partial_t u)^2 r^{d-1}\,dr\le 0
    \end{equation}
    ensures that the blow-up can happen only at $r=0$.  Let us define
    $R(t)$ as the smallest spatial scale involved in the blow-up
    (obviously, $R(t)\to0$ with $t\to T$).  When we approach the
    blow-up time, the solution on the scale $r=\mathcal O(R(t))$ looks
    like $u(r,t)=Q\left(\frac{r}{R(t)}\right)$ for some fixed profile
    $Q$.  This motivates the following definition of a blow-up rate
    \begin{equation}
    \label{eq:blow-up-rate}
    R(t)=\frac{1}{\sup_{r\ge0}\,\lvert \partial_r u(r,t) \rvert}.
    \end{equation}
    By the definition \eqref{eq:blow-up-rate} of the blow-up rate
    $R(t)$, a re-scaled solution $u(r/R(t),t)$ has a bounded gradient
    for all times $t<T$:
    \begin{equation*}
    \sup_{r\ge0}\,\lvert \partial_r u(r/R(t),t) \rvert=1.
    \end{equation*}

\label{sec-1-5}

    The blow-up mechanisms governed by \eqref{eq:u} depend heavily on
    $k$ and $d$ and can be either Type I or Type II.  For
    $k$-corotational maps in dimension $d=2$ van den Berg, Hulshof and
    King \cite{VandenBerg2003} derived formal results for blow-up
    rates.  In particular, for $1$-corotational maps, they conjectured
    that the generic blow-up is of Type II with the blow-up rate
    \begin{align*}
    R(t)\sim \frac{(T-t)}{|\log(T-t)|^2}\qquad \text{as }t\nearrow T.
    \end{align*}
    Recently, this result has been proved by Raphael and Schweyer
    \cite{Raphael2011} by using methods coming from analysis of
    dispersive equations.  For $1$-corotational maps in dimension $2$
    other, non-generic blow-up rates, are also possible
    \cite{Angenent2009}.

    For $1$-corotational maps in dimensions $3\le d\le 6$ Fan
    \cite{Fan1999} used ODE methods to prove the existence of a
    countable family $\{f_n\}_{n=1,2,\dots}$ of self-similar
    solutions for which
    \begin{align*}
    R(t)\sim (T-t)^{\frac{1}{2}}.
    \end{align*}
    Later, Biernat and Bizo\'n \cite{Biernat2011} showed, via
    numerical and analytical methods, that only $f_1$ is linearly
    stable and corresponds to a generic Type I blow-up.  Gastel
    \cite{Gastel2002a} proved that the solution $f_1$ exists also for
    $k$-corotational maps as long as $3\le d<2+k(2+2\sqrt{2})$.  On the
    other hand, there are no results in the literature on dimensions
    $d>2+k(2+2\sqrt{2})$, even for $1$-corotational maps.
\subsection*{Statement of the main result}
\label{sec-1-1}

\label{sec-1-1-1}

    In our paper we use a method of matched asymptotics to construct a
    generic type II solution for $1$-corotational maps in dimensions
    $d\ge7$.  As $t\nearrow T$, the blow-up rates of these solutions
    are asymptotically given by
    \begin{align}
    \label{eq:generic-k1-d7}
    R(t)&\sim\frac{C(T-t)^{\frac{1}{2}}}{-\log(T-t)-\kappa}& \text{for }d&=7 \\
    \label{eq:generic-k1}
    R(t)&\sim\kappa (T-t)^{\frac{1}{2}+\beta_1} & \text{for } d&>7
    \end{align}
    with $\beta_1>0$ defined as
    \begin{align}
    \label{eq:3}
    \beta_1=-\frac{1}{2}+\frac{2}{d-2-\omega},\quad \omega=\sqrt{d^2-8d+8}.
    \end{align}
    For each blow-up rate the constant $\kappa$ represents the
    dependence on initial data, while in \eqref{eq:generic-k1-d7} the
    constant $C$ is a fixed number.  Interestingly, the blow-up rate
    in dimension $d=7$ is, to the leading order, equal to
    \begin{equation*}
    R(t)=\frac{C(T-t)^{\frac{1}{2}}}{-\log(T-t)}(1+\mathcal O(\lvert\log(T-t)\rvert^{-1})),\qquad t\nearrow T
    \end{equation*}
    so, in dimension $7$, the blow-up rate is asymptotically
    independent of initial data.

\label{sec-1-1-2}

    Dimension $d=7$ can be seen as a borderline between type I and
    type II blow-ups.  If one forgets about the underlying geometric
    setup and allows for non-integer values of $d$, then all our
    results remain valid.  For $d$ slightly less than $7$ numerical
    evidence indicates a presence of a generic type I blow-up.  On the
    other hand, when $d$ approaches $7$ from above, $\beta_1$
    continuously drops to zero.  So for $d<7$ we have a type I blow-up
    but for all $d>7$ we have a power-law type II blow-up of the form
    \eqref{eq:generic-k1}.  Naively, one could arrive to a conclusion
    that for $d=7$ we should have a type I blow-up.  Instead, we get a
    type II blow-up \eqref{eq:generic-k1-d7} corresponding to a type I
    blow-up rate with a logarithmic correction.  The transition from
    type I to type II solution at $d=7$ also indicates that the
    self-similar solutions to \eqref{eq:generic-k1-d7} cease to exist
    for $d\ge7$; but analysis of these vanishing self-similar
    solutions is beyond the scope of this paper.

\label{sec-1-1-3}

    In fact, the results for $1$-corotational maps are a special case
    of a more general result for $k$-corotational maps that we
    derive.  For $k$-corotational maps with dimension $d$ and any
    positive integer $N$ satisfying
    \begin{align}
    \label{eq:rate-kn}
    R_N(t)&\sim\kappa(T-t)^{\frac{1}{2}+\beta_N}&\text{for}\qquad
    \begin{split}
    d&>2+k(2+2\sqrt{2})\\
    N&>\frac{1}{4}(d-2-\omega),
    \end{split}\\
    \intertext{}
    \label{eq:rate-kn-neutral}
    R_N(t)&=\frac{C\,(T-t)^{\frac{1}{2}}}{(-\log(T-t)+s_0)^{\frac{1}{\delta}}} & \text{for}\qquad
    \begin{split}
    d&>2+k(2+2\sqrt{2})\\
    N&=\frac{1}{4}(d-2-\omega),
    \end{split}
    \end{align}
    with $\beta_N>0$ defined as
    \begin{align*}
    \beta_N=-\frac{1}{2}+\frac{2N}{d-2-\omega}, \quad \omega=\sqrt{(d-2(k+1))^2-8k^2}
    \end{align*}
    and $\delta>0$ equal to
    \begin{align}
    \label{eq:peculiar-delta}
    \delta=\min(\omega,d-2-\omega).
    \end{align}
    From the dynamical system point of view, each of these solutions
    corresponds to a saddle point with $N-1$ unstable directions.  The
    constants $\kappa$ and $s_0$ depend on initial data, while $C$ is
    a function of $d$ and $k$ only.  This means that asymptotically
    blow-up rate \eqref{eq:rate-kn-neutral} is universal for all
    initial data:
    \begin{equation*}
    R_N(t)=\frac{C\,(T-t)^{\frac{1}{2}}}{(-\log(T-t))^{\frac{1}{\delta}}}(1+\mathcal O(\lvert\log(T-t)\rvert^{-1})),\qquad \text{as } t\nearrow T.
    \end{equation*}

\label{sec-1-1-4}

    To obtain the blow-up rates we employed a technique, called
    matched asymptotics, which allows to construct approximate
    solutions to a differential equation on several spatial
    scales.  The method of matched asymptotics expansions was also
    used to obtain formal type II solutions for the equation
    $\partial_t u=\Delta u+u^p$ in \cite{Herrero1994} (see
    \cite{herrero-velazquez-unpublished} for details); these solutions
    have a similar stability properties as solutions
    \eqref{eq:rate-kn}.  On the other hand, the case of
    $1$-corotational maps in $d=7$ (and \eqref{eq:rate-kn-neutral} in
    general) resembles the solutions found by Herrero and Velázquez
    who used matched asymptotic to derive blow-up rates for chemotaxis
    aggregation in \cite{Herrero1996,Vela2002} and for the problem of
    melting ice balls in \cite{Herrero1997}.

    As in the papers of Herrero and Velázquez, the blow-up rates are
    closely connected to the eigenvalues of a singular self-similar
    solution.  In the case of $k$-equivariant harmonic maps, this
    singular solution is remarkably simple, as it corresponds to a
    singular equatorial map $u(r,t)=\frac{\pi}{2}$.  The eigenvalues
    coming from linearization around the equatorial map ($\lambda_N =
    -\frac{d-2-\omega}{2}+N$ for $N=0,1,2,\dots$, see also
    \eqref{eq:eigen}) relate to the blow-up rate exponents via
    $\beta_N=\frac{\lambda_N}{(d-2-\omega)}$.  The interesting case of
    neutral eigenvalues, $\lambda_N=0$, requires us to include
    non-linear corrections into our analysis and gives rise to the
    logarithmic terms in the blow-up rate
    \eqref{eq:rate-kn-neutral}.  Because there are two ways in which
    the nonlinear term can enter the equation, we have to estimate
    both of them and decide which is the dominant one.  Surprisingly,
    this dominance---and thus the blow-up mechanism---is not set in
    stone but it depends on the dimension, which is reflected by a
    peculiar formula \eqref{eq:peculiar-delta}.

\label{sec-1-1-5}

    The formal solutions constructed in this paper are a first step
    towards the rigorous proof of existence of Type II blow-up for the
    equations of heat flow for $k$-corotational harmonic maps.  The
    solutions presented here will be proved to exist in the upcoming
    paper by the author and Yukihiro Seki \cite{biernat-seki}.  The
    proof bases on topological methods similar to the ones used by
    Herrero and Velazquez in \cite{herrero-velazquez-unpublished}.
\section*{Construction of a blowing up solutions}
\label{sec-2}
\subsection*{Preliminaries}
\label{sec-2-1}

\label{sec-2-1-1}

    To describe blow-up at time $T$ it is convenient to introduce the
    self-similar variables
    \begin{equation}
    \label{eq:self-similar}
    y=\frac{r}{\sqrt{T-t}},\quad s=-\log(T-t),\quad f(y,s)=u(r,t)
    \end{equation}
    in which the original equation (\ref{eq:u}) takes the following
    form
    \begin{equation}
    \label{eq:f}
    \partial_s f=\partial_{yy}f+\left(\frac{d-1}{y}-\frac{y}{2}\right)\partial_{y}f-\frac{k(d+k-2)}{2y^2}\sin(2f)
    \end{equation}
    The boundary condition $u(0,t)=0$ trivially carries over as
    $f(0,s)=0$.

\label{sec-2-1-2}

    Self-similar solutions are stationary points of the above
    equation, if they existm they fully capture the blow-up rate
    (i.e. the solution is regular for all $s$ including
    $s=\infty$).  For $1$-corotational maps a countable family
    $\{f_n\}_{n=1,2,\dots}$ of self-similar solutions was proved to
    exist for $3 \le d \le 6$ by Fan \cite{Fan1999}.  Biernat\&Bizoń
    \cite{Biernat2011} demonstrated that only the first member of the
    family, $f_1$, is linearly stable.  Numerical evidence suggests
    that for $d\ge7$ these solutions are absent and therefore the Type
    I blow-up is no longer possible.  For higher topological degrees
    the only rigorous result on existence of self-similar solutions,
    that authors are aware of, is the one by Gastel \cite{Gastel2002a}
    who proved the existence of the monotone self-similar solution
    $f_1$ in dimensions $d\le 2+k(2+2\sqrt{2})$.  Numerical evidence
    suggests, that for all $k\ge1$ and $d\le 2+k(2+2\sqrt{2})$ there
    exists a countable family of self-similar solutions
    $\{f_n\}_{n=1,2\dots}$.

\label{sec-2-1-3}

    On the other hand, in any dimension $d$ and for any topological
    degree $k$ \eqref{eq:f} there exists a singular stationary
    solution, $f(y,s)=\pi/2$.  This solution is singular because it
    violates the boundary conditions at $y=0$.  Linear stability of
    this solution heavily depends on $d$ and $k$.  For $k=1$ and
    dimension $d\ge7$, $f(y,s)=\pi/2$ is linearly stable (up to a gauge
    mode corresponding to the shift of blow-up time $T$).  For $k\ge2$
    and $d>2+k(2+2\sqrt{2})$, $f(y,s)=\pi/2$ looses some stability and
    becomes a saddle point with a finite number of unstable
    directions.  As we shall see, this solution plays the key role in
    the dynamics of the blow-up.
\subsection*{Boundary layer}
\label{sec-2-2}

\label{sec-2-2-1}

    The singular solution $f(y,s)=\pi/2$ serves as a starting point
    for our construction of a Type II blow-up.  The first step is to
    assume that the constructed solution converges to $\pi/2$.  The
    convergence to $\pi/2$ has to be non-uniform because of the
    boundary condition at the origin $f(0,s)=0$.  The non-uniform
    convergence can be realized by a boundary layer of size
    $\epsilon(s)$ near the origin, where a rapid transition from $f=0$
    to $f=\pi/2$ occurs.  This transition can be described by changing
    variables in \eqref{eq:f} to
    \begin{equation}
    \label{eq:xi}
    \xi=\frac{y}{\epsilon(s)},\quad U(\xi,s)=f(y,s),
    \end{equation}
    where the dependent variable $U$ solves
    \begin{equation}
    \label{eq:U}
    \epsilon^2\partial_s U = \partial_{\xi\xi}U+\left(\frac{d-1}{\xi}+(2\epsilon\dot\epsilon-\epsilon^2)\frac{\xi}{2}\right)\partial_\xi U-\frac{k(d+k-2)}{2\xi^2}\sin(2U),\quad U(0,s)=0.
    \end{equation}

\label{sec-2-2-2}

    We expect convergence to $\pi/2$, so the width of the boundary
    layer must tend to zero with time, hence $\epsilon(s)\to0$ for
    $s\to\infty$.  Additionally, we assume that the derivative of
    $\epsilon$ is bounded by $\epsilon$ for large $s$ i.e. $\dot
    \epsilon(s)=\mathcal O(\epsilon(s))$ as $s\to\infty$. Under these
    assumptions one can drop the quadratic terms in $\epsilon$ and
    $\dot\epsilon$ from equation \eqref{eq:U}.  This leads to a
    solution $U(\xi,s)=U^*(\xi)$, where $U^*(\xi)$ solves an ordinary
    differential equation
    \begin{equation}
    \label{eq:u_stat}
    \frac{d^2U^*}{d\xi^2}+\frac{d-1}{\xi}\frac{dU^*}{d\xi}-\frac{k(d+k-2)}{2\xi^2}\sin(2U^*)=0
    \end{equation}
    with boundary condition $U^*(0)=0$ inherited from
    \eqref{eq:U}.  Any $U^*$ solving \eqref{eq:u_stat} is also a
    stationary point of \eqref{eq:u}, i.e. $U^*$ is a $k$-corotational
    harmonic map.

\label{sec-2-2-3}

    Equation \eqref{eq:u_stat} possesses a scaling symmetry
    $\xi\to\lambda\xi$ (with $\lambda>0$), which implies that any
    $U_\lambda(\xi,s)=U^*(\lambda\xi)$ is also an admissible
    approximate solution to \eqref{eq:U}.  To get rid of this
    ambiguity, we first notice, that any regular solution to
    \eqref{eq:u_stat} behaves like $U^*(\xi)=a\xi^k+\mathcal
    O(\xi^{3k})$ near the origin with some real $a$.  We can fix the
    scaling freedom by setting $a=1$, or equivalently by introducing
    an additional boundary condition
    \begin{equation}
    \label{eq:u-boundary}
    U^*(\xi)=\xi^k+\mathcal O(\xi^{3k}) \qquad \text{as } \xi\to0.
    \end{equation}

\label{sec-2-2-4}

    Equation \eqref{eq:u_stat} simplifies to an autonomous system if
    we use variables $x$ and $v$ defined as $\xi=e^{x}$ and
    $2U^*(\xi)=\pi+v(x)$
    \begin{equation}
    \label{eq:v}
    v''+(d-2)v'+k(d+k-2)\sin(v)=0.
    \end{equation}
    The boundary condition $U^*(\xi)=\xi^k+\mathcal O(\xi^{3k})$
    implies $v(x)=-\pi+e^{kx}+\mathcal O(e^{3kx})$ when $x\to
    -\infty$.  Because \eqref{eq:v} is an autonomous equation we can
    deduce some global properties of $U^*$ by analyzing the phase
    diagram of $\eqref{eq:v}$.

\begin{figure}[htb]
    \centering
    \includegraphics[width=.5\linewidth]{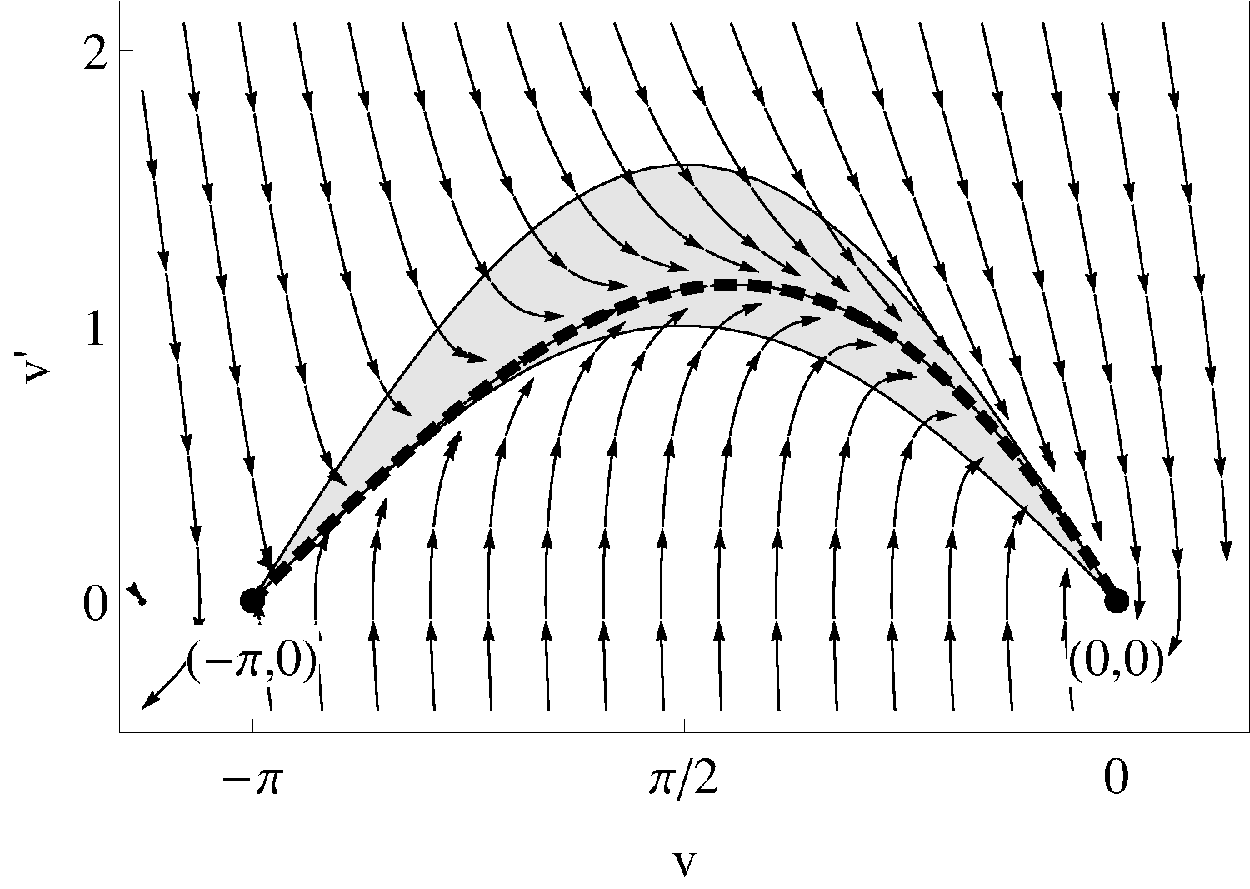}
    \caption{\label{fig:phase}Phase diagram for the equation $0=v''+(d-2)v'+k(d+k-2)\sin(v)$ with $k=1$ and $d=8$.  A solution joining two critical points of the phase diagram is shown as a dashed line.  Additionally, the plot depicts a trapping region $\mathcal S=\{(v,v')\,|\,k\sin(v)\le v'\le \gamma\sin(v)\}$, from which no solution can escape (here $\gamma=\frac{1}{2}(d-2-\sqrt{(d-2(k+1))^2-8k^2})=3-\sqrt{2}$).  The trapping region $\mathcal S$ is used to prove estimates on a depicted solution in Theorem \ref{th:v-asymp}.}
    \end{figure}

    The solution to \eqref{eq:v} subject to these boundary conditions
    has a mechanical interpretation of a motion of a damped pendulum
    with $v$ being the angular position and $x$ corresponding to the
    time.  The boundary conditions demand that the pendulum starts
    inverted, $v=-\pi$, at time $x=-\infty$ and swings out of this
    unstable position.  The damping term forces the pendulum to reach
    the bottom, $v=0$, when $x=\infty$.  In the phase plane spanned by
    $(v,v')$, this trajectory is a heteroclinic orbit starting at the
    saddle point $(-\pi,0)$ and ending at $(0,0)$.  To get the
    asymptotic behavior of $U^*$ at $\xi\to\infty$ it is enough to
    linearize \eqref{eq:v} at the endpoint of the heteroclinic orbit,
    as shown in Figure \ref{fig:phase}.

\label{sec-2-2-5}

    To analyze the asymptotic behavior of $v(x)$ for $x\to\infty$ we
    linearize the equation \eqref{eq:v} at the stationary point
    $(0,0)$.  The eigenvalues of the linearized equation are
    \begin{equation}
    \label{eq:omega}
    \mu_+=-\gamma,\qquad\mu_-=-\gamma-\omega
    \end{equation}
    with constants $\gamma=\frac{1}{2}(d-2-\omega)$ and
    $\omega=\sqrt{(d-2(k+1))^2-8k^2}$.  From the form of the
    eigenvalues $\mu_\pm$ we see that the stationary point $(0,0)$ is
    a stable spiral for $d<d^*:=2+k(2+2\sqrt{2})$ but changes to a
    stable node when $d\ge d^*$.  It follows that the asymptotic
    behavior of $v$, and consequently of $U^*$, can be either
    oscillatory or non oscillatory depending on $d$ for a given
    $k$.  To proceed with our construction, we have to assume the
    latter---non oscillatory---behavior, that is $d\ge d^*$.  For the
    particular case of $1$-corotational maps this condition simplifies
    to $d\ge 7$, if we consider only integer values of $d$.

\label{sec-2-2-6}

    There is one last thing to establish before we can make a claim
    about the asymptotic behavior of $U^*$.  The formula for
    asymptotic behavior of $v$ near $(0,0)$, written explicitly, is
    \begin{equation}
    v(x)=2h_+\cdot e^{x\mu_+}(1+\mathcal O(e^{-2x}))+2h_-\cdot e^{x\mu_-}(1+\mathcal O(e^{-2x}))
    \end{equation}
    (the factor of $2$ is a matter of convenience).  Because
    $\mu_+>\mu_-$, the leading order term should be $2h_+ e^{x\mu_+}$,
    unless $h_+$ is zero, in which case the dominant behavior changes
    to $2h_-e^{x\mu_-}$.  In the appendix (Theorem \ref{th:v-asymp})
    we exclude this possibility by proving that that $h_+$ is
    negative.  We finally conclude that the asymptotic behavior of
    $U^*$ for large $\xi$ is
    \begin{equation}
    \label{eq:U_asymp}
    U^*(\xi)=\frac{\pi}{2}-h\xi^{-\gamma}(1+\mathcal O(\xi^{-2})+\mathcal O(\xi^{-\omega})),
    \end{equation}
    with $h=-h_+>0$ depending only on $d$, and $\gamma>0$ defined as
    $-\mu_+$:
    \begin{equation}
    \label{eq:gamma}
    \gamma=\frac{1}{2}(d-2-\omega),\quad \omega=\sqrt{d^2-8d-8}.
    \end{equation}

\label{sec-2-2-7}

    Let us check where the approximation of $U(\xi,s)$ by $U^*(\xi)$
    is valid.  To arrive at the approximate equation \eqref{eq:u_stat}
    we had to drop the terms containing $\epsilon$ and
    $\dot\epsilon=\mathcal O(\epsilon)$.  The approximation fails if
    one of the dropped terms becomes comparable with the remaining
    terms.  For example we assumed that the remainder term in
    \begin{equation*}
    \left(\frac{d-1}{\xi}+(2\epsilon\dot\epsilon-\epsilon^2)\frac{\xi}{2}\right)=\frac{d-1}{\xi}(1+\mathcal O(\epsilon^2))
    \end{equation*}
    is small.  But this assertion clearly fails for $\xi$ of order
    $1/\epsilon$, so the approximation $U(\xi,s)\approx U^*(\xi)$ can be
    valid only if $\xi\ll 1/\epsilon$ or, by definition \eqref{eq:xi},
    if $y\ll 1$.
\subsection*{Linearization around the singular solution}
\label{sec-2-3}

\label{sec-2-3-1}

    The boundary layer from the previous section resolves a conflict
    between the boundary condition $f(0,s)=0$ and the assumed
    convergence of $f(y,s)$ to $\pi/2$.  In this section, we focus on
    describing the solution to \eqref{eq:f} away from the boundary
    layer, i.e. for $y$ of order $1$.  For such $y$, we expect the
    solution to stay close to $f=\pi/2$, so it is convenient to
    introduce a new variable $\psi$ defined as

    $$
    f(y,s)=\pi/2+\psi(y,s).
    $$

    The new variable $\psi$ solves

    \begin{equation}
    \label{eq:psi}
    \partial_s \psi=-\mathcal A\psi + F(\psi),\qquad F(\psi)=\frac{k(d+k-2)}{2y^2}\left(\sin(2\psi)-2\psi\right)=\mathcal O(\psi^3)
    \end{equation}

    with operator $\mathcal A$ given by

    $$
    -\mathcal A\psi=\frac{1}{\rho}\partial_y\left(\rho\partial_y\psi\right)+\frac{k(d+k-2)}{y^2}\psi,\quad \rho(y)=y^{d-1}e^{-y^2/4}.
    $$

    A natural Hilbert space, arising in the context of operator
    $\mathcal A$ is

    $$
    L^2(\mathbb R_+,\rho\,dy)=\left\{f\in L^2_{loc}(\mathbb R_+)\,\,|\,\,
    \int_0^\infty f(y)^2\rho(y)\,dy < \infty\right\}
    $$

    with a canonical inner product

    \begin{equation}
    \label{eq:dot_product}
    \langle f,g \rangle = \int_0^\infty f(y)g(y)\rho(y)\,dy.
    \end{equation}

    It is routine to check that the operator $\mathcal A$, under the
    assumption $d>2+k(2+2\sqrt{2})$, is self-adjoint in $L^2(\mathbb
    R_+,\rho\,dy)$ with domain $H^1(\mathbb R_+,\rho\,dy)$ --- a
    weighted Sobolev space defined in a canonical way.

\label{sec-2-3-2}

    To find the eigenfunctions of $\mathcal A$ we have to solve an
    ordinary differential equation
    \begin{equation}
    \label{eq:eigen_ode}
    \frac{1}{\rho}\frac{d}{dy}\left(\rho\, \frac{d}{dy}\phi\right)+\frac{k(d+k-2)}{y^2}\phi=-\lambda \phi
    \end{equation}
    with the condition $\phi\in H^1(\mathbb R_+,\rho(y)\,dy)$.  After
    a change of variables $\phi(y)=y^{-\gamma}w(y^2/4)$ and $z=y^2/4$
    (with $\omega$ and $\gamma$ defined in \eqref{eq:gamma}) equation
    \eqref{eq:eigen_ode} becomes
    \begin{equation}
    \label{eq:eigen_z}
    z\frac{d^2w}{dz^2}+\left(1-z+\frac{\omega}{2}\right)\frac{dw}{dz}=-\left(\lambda+\frac{\gamma}{2}\right)w.
    \end{equation}
    with the condition $w\in H^1(\mathbb
    R_+,e^{-z}z^{1+\omega/2}\,dz)$.  Combination of the latter
    condition and the eigenvalue problem \eqref{eq:eigen_z} leads to
    $w(z)=L^{(\omega/2)}_n(z)$ with $\lambda_n+\gamma/2=n$
    ($n=0,1,2,\dots$), where $L^{(\alpha)}_n(z)$ denotes associated
    Laguerre polynomials.  In terms of $\phi$ and $y$ these results
    read
    \begin{equation}
    \label{eq:eigen}
    \phi_n=\mathcal N_n y^{-\gamma} L_n^{(\omega/2)}(y^2/4),\quad \lambda_n=-\gamma/2+n,\quad n=0,1,2,\dots
    \end{equation}
    The normalization constant
    \begin{equation}
    \mathcal N_n=2^{-1-\omega/2}\sqrt{\frac{\Gamma(n+1)}{\Gamma(n+1+\omega/2)}}
    \end{equation}
    assures the orthonormality condition $\langle\phi_n,\phi_m
    \rangle=\delta_{n,m}$.  For completeness we shall add that the
    behavior of $\phi_n$ near the origin is
    \begin{equation}
    \label{eq:phi-origin}
    \phi_n=c_n y^{-\gamma}(1+\mathcal O(y^{-2})),\quad c_n=\frac{2^{-1-\omega/2}}{\Gamma(1+\omega/2)}\sqrt{\frac{\Gamma(1+n+\omega/2)}{\Gamma(1+n)}}.
    \end{equation}

\label{sec-2-3-3}

    Given the orthogonality relation and completeness of $\phi_n$ we
    can represent any solution to \eqref{eq:psi} as the following
    series
    \begin{equation}
    \label{eq:psi_solution}
    \psi(y,s)=\sum_{n=0}^{\infty} a_n(s)\phi_n(y),
    \end{equation}
    In the above expression $a_n(s)$ solve non-linear equations
    \begin{equation}
    \label{eq:an}
    \dot a_n=-\lambda_n a_n+\langle F(\psi),\phi_n\rangle \quad \text{for } n=0,1,2,\dots
    \end{equation}
    with $\dot a_n$ standing for the derivative of $a_n$ with respect
    to $s$ and $F(\psi)$ is defined in \eqref{eq:psi}.  Unfortunately,
    the presence of the non-linear coupling term $\langle
    F(\psi),\phi_n\rangle$ renders \eqref{eq:an} impossible to solve
    in its current form.  In the next section we will make assumptions
    on the form of $\psi$, that will allow us to estimate the non
    linear term.  Consequently, we will be able to produce an
    approximate solution to \eqref{eq:f}.
\subsection*{Construction of a global solution}
\label{sec-2-4}

\label{sec-2-4-1}

    The analysis of the boundary layer solution gives us an
    approximation
    \begin{equation}
    \label{eq:f_inner}
    f(y,s)\approx f_{inn}(y,s)=U^*\left(\frac{y}{\epsilon(s)}\right)\quad \text{for}\quad y\ll 1.
    \end{equation}
    If we take $\epsilon \ll y \ll 1$ we can use the
    asymptotic formula \eqref{eq:U_asymp} for $U^*$ to get
    \begin{equation}
    \label{eq:f_inner_asymp}
    f_{inn}(y,s)= \frac{\pi}{2}-h\epsilon(s)^\gamma y^{-\gamma}
    \end{equation}
    to the leading order.  Because $\epsilon(s)\to0$ with
    $s\to\infty$, the inner solution $f_{inn}(y,s)$ can get
    arbitrarily close to $\pi/2$ for a fixed $y$.  But if $f(y,s)$ is
    close to $\pi/2$ the eigenfunctions of the linear operator
    $\mathcal A$ should work as a good approximation to the solution
    $f(y,s)$, so we write
    \begin{equation}
    \label{eq:f_outer}
    f(y,s)\approx f_{out}(y,s)=\frac{\pi}{2}+\sum_{n=0}^\infty a_n(s)\phi_n(y)\quad \text{for}\quad y\gg \epsilon(s).
    \end{equation}

    Without further assumptions, equations \eqref{eq:an} for the
    coefficients $a_n$ cannot be solved.  To proceed with our
    construction we have to reduce the number of independent degrees
    of freedom; we achieve this by assuming that one coefficient, say
    $a_N$, dominates the others, i.e.
    \begin{equation}
    \label{eq:aN-assump}
    \lvert a_N(s)\rvert\gg \lvert a_n(s)\rvert \qquad \text{for }n\ne N,\quad \text{and}\quad s\to\infty.
    \end{equation}
    By \eqref{eq:aN-assump} the outer solution is dominated by only
    one eigenfunction $\phi_N$ for large $s$
    \begin{equation}
    \label{eq:f_outer}
    f(y,s)\approx f_{out}(y,s)=\frac{\pi}{2}+a_N(s)\phi_N(y)\quad \text{for}\quad y\gg \epsilon(s).
    \end{equation}
    So far, this is the most arbitrary assumption we make, so it is
    critical to ensure that it does not lead to a contradiction at the
    end of the construction.  In one of the following sections we
    verify this assumption and show which conditions on initial data
    does \eqref{eq:aN-assump} require.  This analysis leads to
    conclusions regarding the stability of constructed solutions.

\label{sec-2-4-2}

    Both approximations $f_{inn}$ and $f_{out}$ are compatible in the
    region $\epsilon \ll y\ll 1$ if we impose a relation between
    $a_N(s)$ and $\epsilon(s)$.  Indeed, the outer solution behaves
    like
    \begin{equation}
    \label{eq:f_outer_asymp}
    f_{out}(y,s)=\frac{\pi}{2}+a_N(s)\phi_N(y)=\frac{\pi}{2}+c_N a_N(s)  y^{-\gamma}(1+\mathcal O(y^2)).
    \end{equation}
    (cf. \eqref{eq:phi-origin}) and by comparing
    \eqref{eq:f_outer_asymp} with \eqref{eq:f_inner_asymp} we can
    choose $\epsilon$ such that
    \begin{equation}
    \label{eq:matching}
    c_N a_N(s)=-h\epsilon(s)^{\gamma}.
    \end{equation}
    Equation \eqref{eq:matching} is called the matching condition and
    it serves as a link between the inner solution and the outer
    solution.

\label{sec-2-4-3}

\begin{figure}[htb]
    \centering
    \includegraphics[width=.6\linewidth]{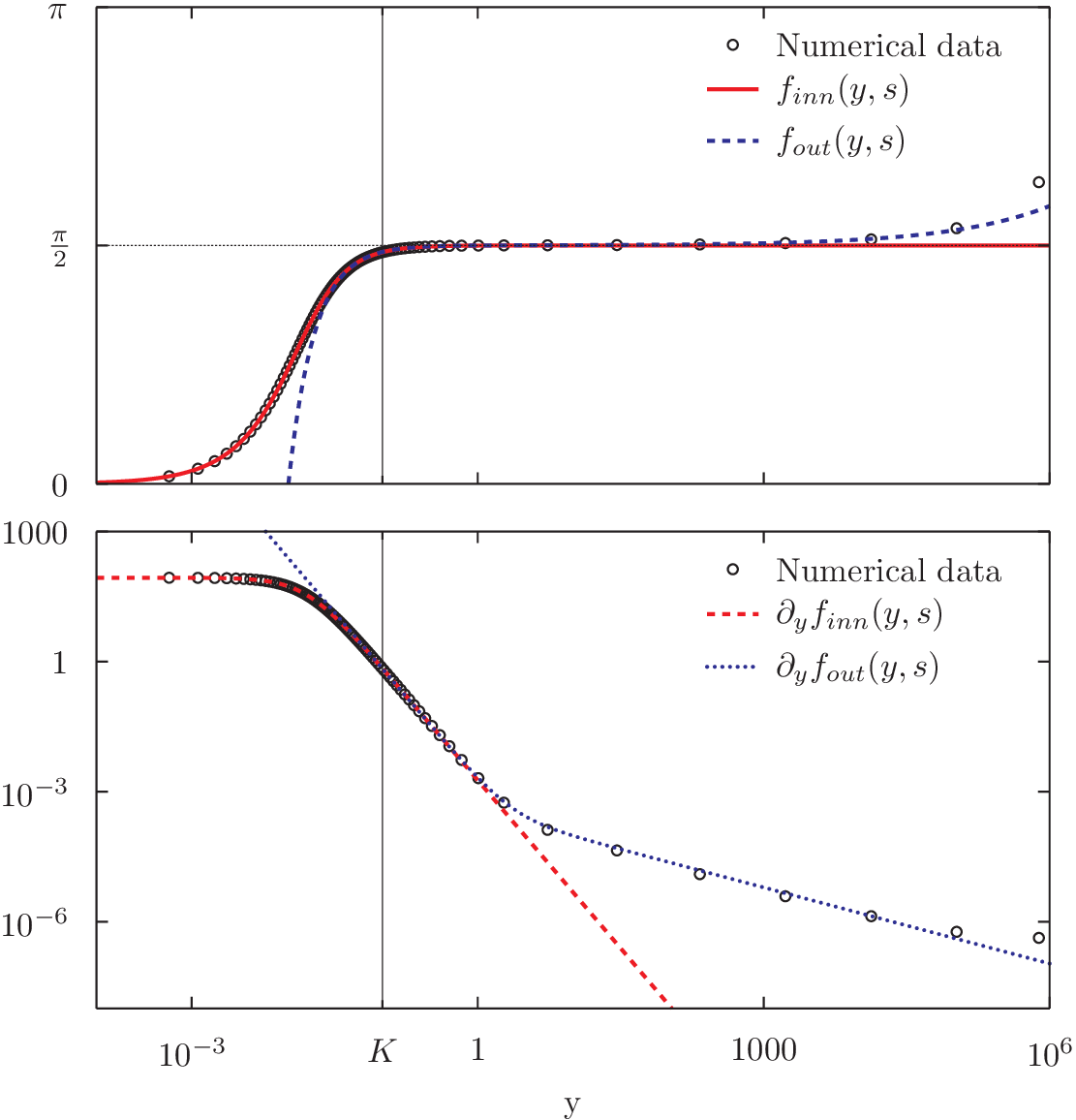}
    \caption{\label{fig:matching}A snapshot at $s=13$ of a numerical solution $f(y,s)$ (in dimension $d=8$) compared to an approximation via inner and outer solutions combined into $f_1(y,s)$ (cf. definition \eqref{eq:f_global}).  The inner solution, $f_{inn}(y,s)$, is a good approximation for $y\ll1$, while the outer solution, $f_{out}(y,s)$, is a good approximation when $f(y,s)$ is close to $\pi/2$.  Both solutions coincide near a point $y=K=10^{-1}$.}
    \end{figure}

    Given solutions \eqref{eq:f_inner} and \eqref{eq:f_outer},
    together with condition \eqref{eq:matching}, we can construct a
    global approximate solution, which is valid for all $y$,
    \begin{equation}
    \label{eq:f_global}
    f_N(y,s)=
    \begin{cases}
    f_\text{inn}(y,s)=U^*\left(\frac{y}{\epsilon(s)}\right)&\text{ for } y\le K\\
    f_\text{out}(y,s)=\pi/2-\frac{h}{c_N}\epsilon(s)^\gamma\phi_{N}(y)&\text{ for } y>K
    \end{cases}
    \end{equation}
    with $K$ chosen so that $\epsilon\ll K\ll 1$
    (e.g. $K=\sqrt{\epsilon}$).  For an example of $f_N$ see Figure
    \ref{fig:matching}.

    At this point, we have an ansatz for a global solution with one
    unknown --- function $\epsilon$.  To get $\epsilon$ we
    have to go back to \eqref{eq:an}, with $n=N$ and
    $a_N(s)=-\epsilon(s)^\gamma\,h/c_N$ and solve
    \begin{equation}
    \label{eq:eps_0}
    \gamma\dot\epsilon=-\lambda_N\epsilon-\,\frac{c_N}{h}\,\epsilon^{1-\gamma}\langle F(\psi),\phi_N\rangle\qquad\psi=f_N(y,s)-\pi/2.
    \end{equation}

\label{sec-2-4-4}

    The remaining question is in what way does the non-linear term
    $\langle F(\psi),\phi_N\rangle$ enter the equation?  To answer
    this question we have to split $\langle F(\psi),\phi_N\rangle$
    into contributions from inner and outer solutions.  However, these
    computations are too technical for this section and would break
    the flow of the argument.  Instead, we enclose the derivation in
    the next section and present the resulting formula here
    \begin{equation}
    \label{eq:nonlinear-estimate}
    \langle F(\psi),\phi_N\rangle = D_N\,\epsilon^{\gamma+\delta},\quad\delta=\min(\omega,2\gamma)>0,\quad  D_N>0.
    \end{equation}
    Combination of the estimate \eqref{eq:nonlinear-estimate} and the
    equation \eqref{eq:eps_0} yields the following equation for
    $\epsilon$
    \begin{equation}
    \label{eq:eps_final}
    \gamma\dot \epsilon=-\lambda_N \epsilon - \frac{D_N c_N}{h}\epsilon^{1+\delta}.
    \end{equation}
    We can immediately discard negative eigenvalues $\lambda_N$, as
    they lead to $\epsilon$ which does not tend to zero; such
    $\epsilon$ violates our previous assumptions about the boundary
    layer.

    The only viable solutions are those with $\lambda_N\ge0$, which
    leads to two further cases.  When $\lambda_N>0$ the non-linear
    term is of higher order and can be discarded for $s$ large
    enough leading to
    \begin{equation}
    \label{eq:epsgt0}
    \epsilon(s)=\epsilon_0\,e^{-\frac{\lambda_N}{\gamma}\,s} \qquad \text{for }\lambda_N>0
    \end{equation}
    with $\epsilon(0)$ depending on initial data.  On the other hand,
    when $\lambda_N=0$ the non-linear term becomes the leading order
    term resulting in
    \begin{equation}
    \label{eq:eps_neutral}
    \epsilon(s)=\frac{C_N}{(s-s_0)^{\frac{1}{\delta}}},\qquad C_N=\left(\frac{h\gamma}{c_N\, D_N\, \delta}\right)^{\frac{1}{\delta}}\qquad \text{for }\lambda_N=0.
    \end{equation}

    We can now relate the blow-up rate $R(t)$ with $\epsilon$ via
    \begin{equation}
    \label{eq:R-epsilon}
    R(t)=\frac{1}{\sup_{r\ge0}\lvert\partial_r u(r,t)\rvert}=C_s\sqrt{T-t}\,\epsilon(s),\qquad C_s=\frac{1}{\sup_{\xi\ge0}\lvert \frac{d U^*}{d\xi}(\xi)\rvert}.
    \end{equation}
    If we combine \eqref{eq:R-epsilon}, \eqref{eq:epsgt0} and
    \eqref{eq:eps_neutral}, and solve the conditions $\lambda_N>0$ and
    $\lambda_N=0$ for $N$ we get the following blow-up rates
    \begin{align}
    \label{eq:rate-kn}
    R_N(t)&=C_s\epsilon_0(T-t)^{\frac{1}{2}+\beta_N} & \text{for}\qquad
    \begin{split}
    d&>2+k(2+2\sqrt{2})\\
    N&>\frac{1}{4}(d-2-\omega)
    \end{split}\\
    \intertext{}
    \label{eq:rate-kn-neutral}
    R_N(t)&=\frac{C_s\,C_N(T-t)^{\frac{1}{2}}}{(-\log(T-t)-s_0)^{\frac{1}{\delta}}} &\text{for}\qquad
    \begin{split}
    d&>2+k(2+2\sqrt{2})\\
    N&=\frac{1}{4}(d-2-\omega),
    \end{split}
    \end{align}
    with $\beta_N>0$
    \begin{align}
    \label{eq:beta-k2}
    \beta_N=-\frac{1}{2}+\frac{2N}{d-2-\omega}, \quad \omega=\sqrt{(d-2(k+1))^2-8k^2}
    \end{align}
    and $\delta>0$ being equal to
    \begin{align}
    \label{eq:delta}
    \delta=\min(\omega,d-2-\omega).
    \end{align}
\subsection*{Approximation of the coupling term}
\label{sec-2-5}

   According to the assumed form of the global solution $f_N(y,s)$ we
   can approximate the solution $\psi$ in the intervals $y\le K$ and
   $y>K$ separately.  Therefore, we split the integral $\langle
   F(\psi),\phi_{n}\rangle$ into
   \begin{equation*}
   \langle F(\psi),\phi_{n}\rangle=\left(\int_{0}^{K}+\int_{K}^{\infty}\right)F(\psi)\phi_{n}(y)y^{d-1}e^{-y^2/4}\,dy= I_{inn}+ I_{out}.
   \end{equation*}
   We compute the two integrals $I_{inn}$ and $I_{out}$ and compare
   them to see which one gives the leading order contribution.  Our
   analysis leads to two qualitatively different approximations of the
   non-linear term
   \begin{equation*}
   F(\psi)=\frac{k(d+k-2)}{2y^2}(\sin(2\psi)-2\psi)
   \end{equation*}
   depending on the choice of $d$ and $k$.

\label{sec-2-5-1}

    The first integral, $I_{inn}$, contains the contribution from the
    inner layer, where $\psi\approx f_{inn}-\pi/2$, so by \eqref{eq:f_inner} we
    can approximate $F(\psi)$ as
    \begin{align*}
    F(\psi)=F(f_{inn}(y,s)-\pi/2)=F(U^*(y/\epsilon)-\pi/2)=y^{-2}g(y/\epsilon),
    \end{align*}
    for brevity we use a notation
    \begin{equation*}
    g(\xi)=\frac{k(d+k-2)}{2}(\sin(2U^*(\xi)-\pi)-(2U^*(\xi)-\pi)).
    \end{equation*}
    When $y<K\ll1$ we can replace the eigenfunction and the weight
    $\phi_{n}(y)\,y^{d-1}e^{-y^{2}/4}$ with its leading order term
    $\phi_{n}(y)\,y^{d-1}e^{-y^{2}/4}=c_{n}y^{-\gamma+d-1}(1+\mathcal
    O(y^{2}))$.  We finally arrive at a simplified version of the
    integral $I_{inn}$
    \begin{align}
    \label{eq:I-inn}
    I_{inn}\approx c_{n}\int_0^K
    g(y/\epsilon)y^{d-3-\gamma}\,dy=c_{n}\epsilon^{d-2-\gamma}\int_0^{K/\epsilon}g(\xi)\xi^{d-3-\gamma}\,d\xi.
    \end{align}
    The upper bound $K/\epsilon$ in \eqref{eq:I-inn} tends to infinity
    as $s\to\infty$, so it is reasonable to check whether the
    integrand is divergent or convergent as $\xi\to\infty$.  To this
    end we have to compute the asymptotic behavior of $g(\xi)$ at
    infinity.  This can be done by using the asymptotic of $U^*$, as
    given by \eqref{eq:U_asymp}
    \begin{align*}
    g(\xi)\approx-\frac{2k(d+k-2)}{3}(U^*(\xi)-\pi/2)^{3}\approx\frac{2k(d+k-2)h^{3}}{3}\xi^{-3\gamma}\quad \text{as }\xi\to\infty
    \end{align*}
    The leading order of the integrand is thus $\xi^{d-3-4\gamma}$.  By
    definitions \eqref{eq:gamma} of $\gamma$ and $\omega$ there holds
    \begin{equation}
    \label{eq:gamma-relation}
    d-2-\gamma=\gamma+\omega,
    \end{equation}
    so the leading order term can be written as
    $\xi^{d-3-4\gamma}=\xi^{\omega-2\gamma-1}$.

    We have to consider two cases, because the integral
    \eqref{eq:I-inn} can be divergent or convergent for large
    $(K/\epsilon)$ depending on the sign of $\omega-2\gamma$.  If
    $\omega < 2\gamma$, then the integral converges so, by taking the
    limit $K/\epsilon\to\infty$, we get
    \begin{align*}
    I_{inn}=\,&c_n\epsilon^{d-2-\gamma}\int_0^{K/\epsilon}g(\xi)\xi^{d-3-\gamma}\,d\xi\\
    =\,&c_n\epsilon^{\gamma+\omega}\int_0^{\infty}g(\xi)\xi^{d-3-\gamma}\,d\xi.
    \end{align*}
    But when $\omega > 2\gamma$ the integral diverges as
    $(K/\epsilon)^{\omega-2\gamma}$, so we can replace the integral
    with its rate of divergence, in which case the lowest order
    approximation is
    \begin{align*}
    I_{inn}\approx
    \frac{2k(d+k-2)h^{3}c_{n}}{3}\epsilon^{\gamma+\omega}\left(\frac{K}{\epsilon}\right)^{\omega-2\gamma}=\frac{2k(d+k-2)h^{3}c_{n}}{3}\epsilon^{3\gamma}K^{\omega-2\gamma}.
    \end{align*}

\label{sec-2-5-2}

    We have to consider two similar cases when dealing with
    $I_{out}$.  For $I_{out}$, $\psi$ is dominated by its
    approximation via a single eigenfunction
    $\psi=-\frac{h}{c_N}\epsilon^\gamma\phi_N$, which, together with
    $y>K$, results in $\lvert\psi\rvert\ll 1$ near the origin.  So, as
    the first step to the approximation of $I_{out}$ we expand $F$ in
    a Taylor series around $\psi=0$
    \begin{align*}
    F(\psi)=F\left(-\frac{h}{c_N}\epsilon^\gamma\phi_N\right)\approx-\frac{k(d+k-2)}{y^2}\cdot\frac{3}{2}\left(-\frac{h}{c_N}\epsilon^\gamma\phi_N\right)^3=\frac{2k(d+k-2)h^3}{3y^{2}c_N^3}\epsilon^{-3\gamma}\phi_{N}^{3}.
    \end{align*}
    If we use the Taylor expansion in $I_{out}$ we obtain
    \begin{align*}
    I_{out}=\frac{2k(d+k-2)h^3}{3c_N^3}\epsilon^{-3\gamma}\int_{K}^{\infty}\phi_{N}(y)^{3}\phi_{n}(y)y^{d-3}e^{-\frac{y^{2}}{4}}\,dy
    \end{align*}
    which can be either divergent or convergent for small $K$.  Near
    the origin ($y\to0$) the leading order behavior of the integrand
    is
    \begin{align*}
    \phi_{N}(y)^{3}\phi_{n}(y)y^{d-3}e^{-\frac{y^{2}}{4}}= (c_N)^3 c_n\, y^{d-3-4\gamma}(1+\mathcal O(y^{2}))=(c_N)^3 c_n\, y^{\omega-2\gamma-1}(1+\mathcal O(y^{2})).
    \end{align*}
    It is clear, that for $\omega<2\gamma$ the integral is finite and
    we can take the limit $K\to0$, while for $\omega>2\gamma$ the
    integral is divergent and behaves like $(c_N)^3c_n
    K^{\omega-2\gamma}$.  These two cases can be expressed as
    \begin{align*}
    I_{out}&\approx\frac{2}{3}k(d+k-2)h^3c_n\,\epsilon^{-3\gamma}\,K^{\omega-2\gamma}\quad &\text{for }\omega<2\gamma,
    \intertext{and}
    \begin{split}
    I_{out}&\approx\, \epsilon^{-3\gamma}\,\frac{2k(d+k-2)h^3}{3c_N^3}\, \int_{0}^{\infty}\phi_{N}(y)^{3} \phi_{n}(y) y^{d-3} e^{-\frac{y^{2}}{4}}\,dy\\
    &=\epsilon^{-3\gamma}\,T_n
    \end{split}
    \quad &\text{for }\omega > 2\gamma.
    \end{align*}

\label{sec-2-5-3}

    We are now in a position to compare the contributions from
    $I_{inn}$ and $I_{out}$
    \begin{align*}
    I_{inn}&\propto \epsilon^{\gamma+\omega},                        & I_{out}&\propto \epsilon^{\gamma+\omega}\,\left(\frac{\epsilon}{K}\right)^{2\gamma-\omega},\quad & \text{for } \omega<2\gamma,\\
    I_{inn}&\propto \epsilon^{3\gamma}K^{\omega-2\gamma}, & I_{out}&\propto \epsilon^{3\gamma},\quad                                                                 & \text{for } \omega>2\gamma.
    \end{align*}
    For sufficiently large times $I_{inn}$ dominates over $I_{out}$
    when $\omega<2\gamma$ because the term
    $(\epsilon/K)^{2\gamma-\omega}$ tends to zero.  On the other hand,
    when $\omega>2\gamma$ it is the other way around and $I_{out}$
    dominates over $I_{inn}$ due to $K^{\omega-2\gamma}\to 0$.  These
    two cases can be written in a unified way as
    \begin{align}
    \label{eq:non-linear}
    \langle F(\psi),\phi_{n}\rangle\approx D_{n}\epsilon^{\gamma+\delta},\quad
    \delta=\min(2\gamma,\omega).
    \end{align}
    with a constant
    \begin{equation*}
    D_n=\begin{cases}
    c_n \int_0^{\infty}g(\xi)\xi^{d-3-\gamma}\,d\xi  & \text{for } \omega<2\gamma\\
    \frac{2(d-1)h^3}{3c_N^3}\, \int_{0}^{\infty}\phi_{N}(y)^{3} \phi_{n}(y) y^{d-3} e^{-\frac{y^{2}}{4}}\,dy
      & \text{for } \omega>2\gamma\,.
    \end{cases}
    \end{equation*}

    We intentionally avoided the case $\omega=2\gamma$, for which both
    integrals diverge logarithmically.  This happens only for a non
    integer dimension $d=\frac{2}{3} \left(7+2 \sqrt{7}\right)\approx
    8.194\dots$, which we can exclude as incompatible with underlying
    geometric setting of the heat flow for harmonic maps.

    One possible interpretation of this phenomenon-is a change of the
    way we should approximate the non-linear term $F(\psi)$ before the
    projection onto $\phi_n$.  For example, when $\omega>2\gamma$, we
    can safely replace $F(\psi)$ with its Taylor expansion near
    $\psi=0$, i.e. $F(\psi)=\frac{2k(d+k-2)}{3y^2}\psi^3$.  Projecting
    $F(\psi)$ back to $\phi_n$ gives negligible contribution from
    $I_{inn}$ and significantly larger contribution from
    $I_{out}$.  At the same time, the value of $I_{out}$ is
    proportional to the third power of amplitude of $\psi\propto a_N$:
    $I_{out}\propto \epsilon^{3\gamma}\propto a_N^3$.

    As for the other case, $\omega<2\gamma$, the contribution from the
    Taylor expansion is subdominant.  Instead, a very small region
    $y<K$, of a diminishing size, governs the leading order behavior
    of non-linear term $F(\psi)$.  We can replicate this effect by
    approximating $F(\psi)$ with a Dirac delta: $F(\psi)=G
    \epsilon^{\gamma+\omega}\delta(y)$.  Indeed, to the leading order
    we get the same values for projections:
    \begin{equation*}
    \langle F(\psi),\phi_n\rangle=G\epsilon^{\gamma+\omega}\langle \delta,\phi_n\rangle=G c_n\epsilon^{\gamma+\omega}.
    \end{equation*}
    In fact, replacing the non-linear term $F(\psi)$ with a Dirac
    delta is the starting point to several derivations of Type II
    solutions\cite{Herrero1997,Herrero1996}.  On the other hand, the
    Taylor expansion rarely shows up in derivations of the blow-up
    rate.

    To verify whether $\langle F(\psi),\phi_N\rangle$ is positive
    (which is required for solution \eqref{eq:eps_neutral}) it
    suffices to show that $D_N>0$.  The first case, when
    $D_N=c_N\int_0^{\infty}g(\xi)\xi^{d-3-\gamma}\,d\xi$ follows from
    the properties of the bounding region used in Theorem
    \ref{th:v-asymp}, which guarantees that $0\le U^*(\xi)<\pi/2$,
    hence $g(\xi)>0$; combined with $c_n>0$ for every $n\ge0$ we get
    the result.  In the second case the result follows from the sign
    of the integrand in
    \begin{equation*}
    \frac{2k(d+k-2)h^3}{3c_N^3}\, \epsilon^{-3\gamma} \int_{0}^{\infty}\phi_{N}(y)^{4}y^{d-3} e^{-\frac{y^{2}}{4}}\,dy>0.
    \end{equation*}
    and from $c_N>0$.
\subsection*{Note on stability of type II solutions}
\label{sec-2-6}

   In this section we address two concerns that arose earlier in the
   text.  The first one is an ex post validation of our assumption
   \eqref{eq:aN-assump} about the dominance of $a_N$ over other
   coefficients $a_n$.  The other issue is the stability of $f_N$.  It
   appears that $f_N$ is unstable, because there is always a negative
   eigenvalue $\lambda_0=-\gamma/2$. To obtain any of the constructed
   solutions we will have to suppress this instability by fine tuning
   of initial data.

\label{sec-2-6-1}

    With an estimate on the non-linear term $\langle
    F(\psi),\phi_n\rangle$, we can actually solve equations
    \eqref{eq:an} for $a_n$.  By plugging
    \eqref{eq:nonlinear-estimate} into \eqref{eq:an} we get linear
    nonhomogeneous equations
    \begin{equation*}
    \dot a_n=-\lambda_n a_n+D_n\epsilon^{\gamma+\delta},\qquad n\ne N
    \end{equation*}
    which can be explicitly solved by
    \begin{equation}
    \label{eq:an-sol-gen}
    a_n(s)=a_n(0)e^{-\lambda_n s}+D_n\int_0^s\epsilon(q)^{\gamma+\delta}e^{-\lambda_n(s-q)}dq.
    \end{equation}
    The free parameters $a_n(0)$ are connected to initial
    data via
    \begin{equation}
    \label{eq:an-initial}
    a_n(0)=\langle\psi,\phi_n \rangle|_{s=0}.
    \end{equation}

\label{sec-2-6-2}

    Let us start with the coefficients in front of higher
    eigenfunctions, i.e. $n>N$.  It is enough to study the limit
    \begin{equation*}
    \lim_{s\to\infty}\frac{a_n(s)}{a_N(s)}=-\frac{c_N}{h}\lim_{s\to\infty}\frac{a_n(0)+D_n\int_0^s\epsilon(q)^{\gamma+\delta}e^{\lambda_n\,q}\,dq}{e^{\lambda_n s}\epsilon(s)^\gamma}
    \end{equation*}
    The denominator diverges to infinity, while the numerator either
    diverges to $\pm\infty$ or converges to a constant.  In the latter
    case the limit is $0$, and we are done.  If the former is true, we
    apply l'Hôpital's rule to get
    \begin{equation*}
    \lim_{s\to\infty}\frac{a_n(s)}{a_N(s)}=-\frac{c_N\,D_n}{h}\lim_{s\to\infty}\frac{\epsilon(s)^{\delta}}{\gamma\dot\epsilon(s)/\epsilon(s)+\lambda_n}\underset{\text{Eq. }\eqref{eq:eps_final}}{=}-\frac{c_N\,D_n}{h}\lim_{s\to\infty}\frac{\epsilon(s)^{\delta}}{(\lambda_n-\lambda_N)-D_N\,c_N\,\epsilon(s)^\delta/h}=0.
    \end{equation*}
    Hence, without any assumptions on $a_n(0)$ we have $\lvert
    a_N(s)\rvert\gg \lvert a_n(s)\rvert$ for $n>N$.

\label{sec-2-6-3}

    For $n<N$, let us rewrite \eqref{eq:an-sol-gen} as
    \begin{equation}
    \label{eq:an-sol-nong}
    a_n(s)=\left(a_n(0)+D_n\int_0^\infty\epsilon(q)^{\gamma+\delta}e^{\lambda_n q}dq\right)e^{-\lambda_n s}-D_n\int_s^\infty\epsilon(q)^{\gamma+\delta}e^{-\lambda_n(s-q)}dq.
    \end{equation}
    With elementary calculations and knowledge of $\epsilon$ one can
    show that the integrals in \eqref{eq:an-sol-nong} converge if
    $n<N$.  The second term in \eqref{eq:an-sol-nong} is actually much
    smaller than $a_N(s)$.  This is evident when we apply l'Hôpital's
    rule to the limit
    \begin{align*}
    \lim_{s\to\infty}\frac{\int_s^\infty\epsilon(q)^{\gamma+\delta}e^{-\lambda_n(s-q)}dq}{a_N(s)}
    &=\lim_{s\to\infty}\frac{\int_s^\infty\epsilon(q)^{\gamma+\delta}e^{\lambda_n q}dq}{e^{\lambda_n s}a_N(s)}\overset{\kern2pt H}{=}\lim_{s\to\infty}\frac{-\epsilon(s)^{\gamma+\delta}}{\dot a_N(s)+\lambda_n a_N(s)}.     \intertext{We continue with the help of matching condition $c_Na_N(s)=-h\epsilon(s)^\gamma$ and equation \eqref{eq:eps_final} for $\epsilon$ to obtain}
    &=-\frac{c_N}{h}\lim_{s\to\infty}\frac{\epsilon(s)^{\delta}}{(\lambda_n-\lambda_N)-D_N c_N\epsilon(s)^{\delta}/h}=0.
    \end{align*}
    In a similar way we can check that for $n<N$ the first term,
    containing $e^{-\lambda_n s}$, is actually much larger than
    $a_N(s)$.  So if we want $\lvert a_N(s)\rvert \gg \lvert
    a_n(s)\rvert$ to hold, the coefficient in front of $e^{-\lambda_n
    s}$ in \eqref{eq:an-sol-nong} has to be zero.  This can be
    accomplished by selecting particular initial data for which
    \begin{equation}
    \label{eq:an-requirement}
    \langle \psi|_{s=0},\phi_n\rangle = a_n(0)=-D_n\int_0^\infty\epsilon(q)^{\gamma+\delta}e^{\lambda_n q}dq\qquad 0\le n<N.
    \end{equation}
    If the initial data, $\psi|_{s=0}$, fulfills the condition
    \eqref{eq:an-requirement} the assumption $\lvert a_N(s)\rvert \gg
    \lvert a_n(s)\rvert$ does not lead to a contradiction.

\label{sec-2-6-4}

    The condition \eqref{eq:an-requirement} for tthe solution $f_N$
    imposes $N$ constraints on the initial data.  Each constraint
    corresponds to one unstable direction along which our solution can
    diverge from the ansatz $f_N$.  There is, however, one free
    parameter---the blow-up time $T$---that we can use to change the
    values of coefficients $a_n(0)$.  Any small change $T\to T+\eta$
    in blow-up time results in a small change of self-similar
    coordinates \eqref{eq:self-similar} $y\to y-\frac{1}{2}\eta e^s
    y+\mathcal O(e^{2s}\eta^2)$ and $s\to s-\eta e^s+\mathcal
    O(e^{2s}\eta^2)$.  This change in self-similar coordinates affects
    the initial data $\psi|_{s=0}$ so the coefficients $a_n(0)$ also
    change.  In particular, the zeroth coefficient becomes
    \begin{align*}
    a_0(0)\,\to\, a_0(0)-\eta\left\langle\partial_s\psi+\frac{y}{2}\partial_y\psi\,,\,\phi_0\right\rangle_{s=0}+\mathcal O(\eta^2).
    \end{align*}
    It should be possible to choose a blow-up time $T$ in such way,
    that the new $a_0(0)$ fulfills the condition
    \eqref{eq:an-requirement}.  This mechanism removes one of the
    constraints on initial data so $f_N$ has effectively $N-1$
    unstable directions.
\subsection*{Discussion of the resultsy}
\label{sec-2-7}

   In the previous section we analyzed the stability of $f_N$
   concluding that the solution $f_N$ has $N-1$ unstable
   directions.  On the other hand, $N$ is constrained by the condition
   $\lambda_N\ge0$, or equivalently,
   \begin{equation}
   \label{eq:N-condition}
   N\ge \frac{1}{4}(d-2-\omega)=\frac{1}{4}(d-2-\sqrt{(d-2(k+1))^2-8k^2})
   \end{equation}
   with $d>2+k(2+2\sqrt{2})$.  The right hand side of the inequality
   \eqref{eq:N-condition} depends on $k$ and $d$ and puts a lower
   bound on the possible $N$.  In turn, the lower bound on $N$ induces
   a condition on the existence of stable $f_N$ for a given $k$.  If
   we take arbitrary $k\ge1$ and $d>2+k(2+2\sqrt{2})$ we can derive a
   lower bound on $N$
   \begin{equation}
   \label{eq:N-stability}
   N\ge \frac{1}{4}(d-2-\omega)=\frac{1}{4}(d-2-\sqrt{(d-2(k+1))^2-8k^2})>\frac{k}{2},
   \end{equation}
   so the instability of solutions $f_N$ increases with topological
   degree $k$.

\label{sec-2-7-1}

    Only a solution with $N=1$ can be stable, so from the bound
    $N>\frac{k}{2}$ we infer that for $k\ge 2$ there are no stable
    solutions $f_N$.  Still, the solutions to \eqref{eq:u} are
    guaranteed to blow up for large class of initial data.  Numerical
    evidence suggests, that the generic blow-up is self-similar, so
    there has to exist at least one self-similar solution to
    \eqref{eq:u} for $k\ge 2$.  However, to the authors knowledge, in
    the literature there are no rigorous results concerning these
    solutions.

\label{sec-2-7-2}

    Because \eqref{eq:N-stability} is only a lower bound, one can ask
    if there are any examples of stable solutions.  Stable solutions
    could exist only for $1$-corotational maps, so let us assume that
    $k=1$.  The lower bound on $d$ then becomes
    $d>2+k(2+2\sqrt{2})=4+2\sqrt{2}\approx 6.828\dots$, but because
    dimension $d$ is an integer we arrive at $d\ge7$.  The first
    eigenvalue $\lambda_1=-\frac{\gamma}{2}+1$, which corresponds to
    the stable solution, is actually positive for all $d\ge7$ so in
    this case there exists a stable solution $f_1$.  In fact,
    numerical evidence suggests, that $f_1$ corresponds to a generic
    blow-up in dimensions $d\ge 7$.

\label{sec-2-7-3}

    Existence of a generic type I solution in a form of $f_1$ can be
    confirmed numerically, although solutions with finite-time
    singularities present several conceptual difficulties when solved
    on a computer.  The most significant problem comes from the
    spatial resolution needed to resolve the shrinking scale of the
    boundary layer.  We overcome this difficulty by employing a well
    established numerical method called a moving mesh, in which a
    constant number of mesh points is distributed dynamically to
    satisfy demands for high mesh density near the singularity, and
    outside of it.  In particular, we modified \cite{Biernat-MOVCOL}
    an existing implementation \cite{Huang-movcol} of moving mesh
    algorithm MOVCOL \cite{Russell2007}.  For an in-depth description
    of an application of MOVCOL to solutions with finite time
    singularity we refer the reader to a paper on a type II blow-up
    for chemotaxis aggregation by Budd et al.\cite{Budd2005}.

\label{sec-2-7-4}

\begin{figure}[htb]
    \centering
    \includegraphics[width=.5\linewidth]{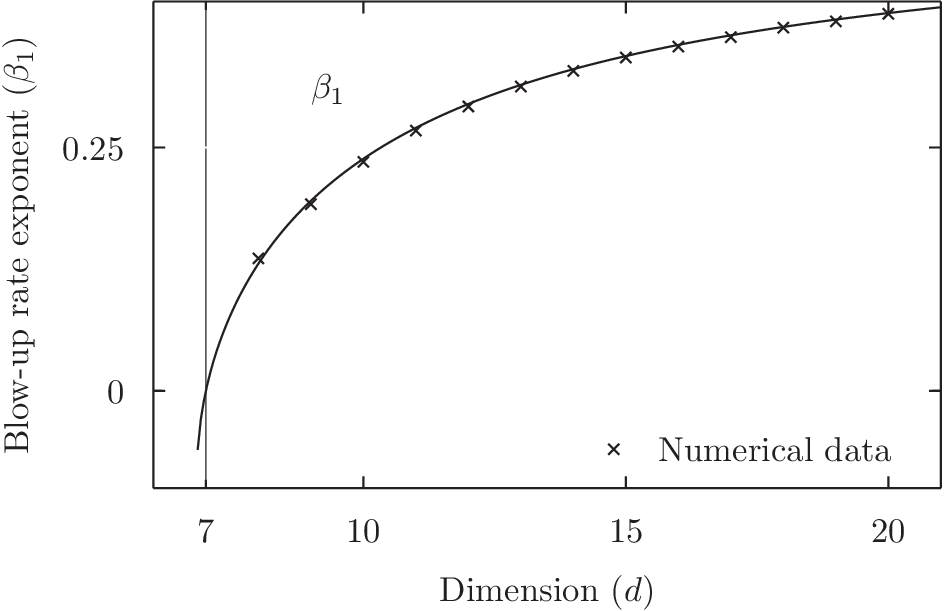}
    \caption{\label{fig:beta1}The predicted blow-up rate for $1$-corotational maps is $R(t)=(T-t)^{\frac{1}{2}+\beta_{1}}$ with $\beta_1=-\frac{1}{2}+\frac{2}{d-2-\omega}$ and $\omega=\sqrt{d^2-8d+8}$.  The figure depicts the comparison between the predicted value of $\beta_1$ and $\beta_1$ obtained from numerical experiment via a relation $\frac{\partial_{tr} u(0,t)}{\partial_r u(0,t)}\to \left(\frac{1}{2}+\beta_1\right)$ with $t\nearrow T$.  In each case the initial data was $u(r,0)=r$.}
    \end{figure}

    For $d\ge 8$ the generic blow-up rate is given by
    \begin{equation}
    \label{eq:d8-num-ratio}
    R(t)=C_s \epsilon_0 (T-t)^{\frac{1}{2}+\beta_1},\qquad \beta_1=-\frac{1}{2}+\frac{2}{d-2-\omega},\qquad \omega=\sqrt{d^2-8d+8}.
    \end{equation}
    By definition \eqref{eq:blow-up-rate} $R(t)$ is inversely
    proportional to $\sup_{r\ge0}\lvert\partial_r u(r,t)\rvert$, which
    can be easily obtained from numerical experiments.  In fact, for
    $k=1$ the supremum is always attained at the point $r=0$, so we
    can replace $\sup_{r\ge0}\lvert\partial_r u(r,t)\rvert$ with
    $\lvert\partial_r u(0,t)\rvert$.  To verify the blow-up rate we
    study the ratio
    \begin{equation*}
    \frac{\partial_{tr} u(0,t)}{\partial_r u(0,t)}=-\frac{R'(t)}{R(t)},
    \end{equation*}
    which in $d\ge8$ should tend to
    \begin{equation*}
    \frac{\partial_{tr} u(0,t)}{\partial_r u(0,t)}\to \left(\frac{1}{2}+\beta_1\right)\qquad \text{as } t\nearrow T.
    \end{equation*}
    We compare $\beta_1$ obtained from numerical experiments with its
    theoretical value in Figure \ref{fig:beta1}.  An additional test
    compares the shape of a numerical solution near the origin with
    the shape of the function $f_1$ with its respective inner and
    outer solutions as in Figure \ref{fig:matching}.  This plot
    captures a solution at time $T-t\approx10^{-5.5}$.

\label{sec-2-7-5}

\begin{figure}[htb]
    \centering
    \includegraphics[width=.7\linewidth]{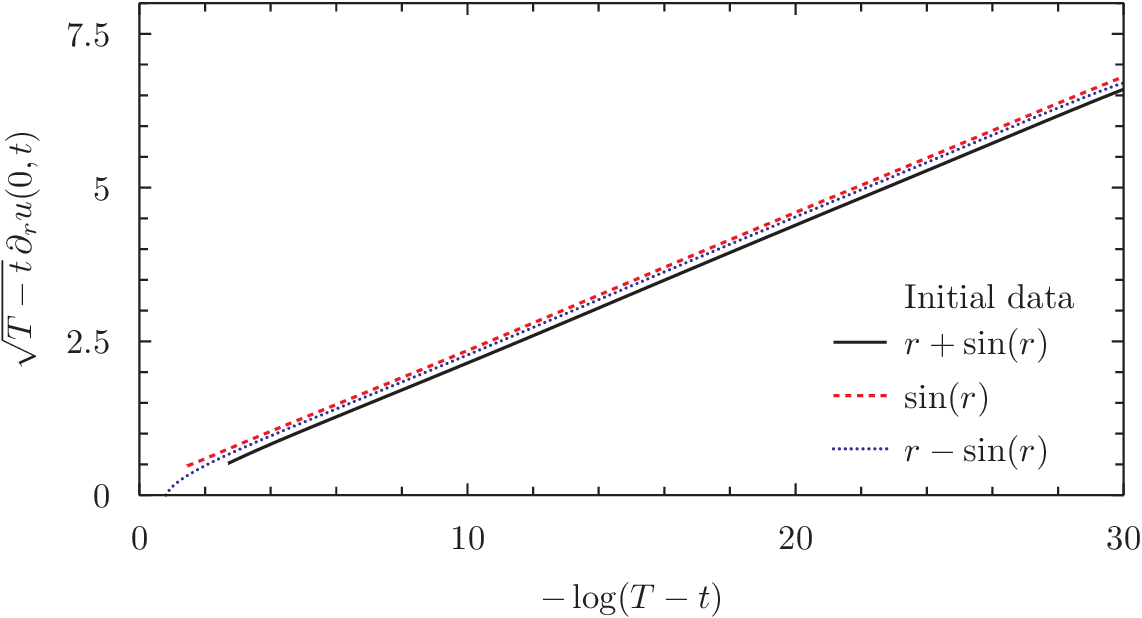}
    \caption{\label{fig:d7}In dimension $d=7$, for a generic blow-up, the rate of blow-up is $R(t)=\frac{C(-\log(T-t)-s_0}{\sqrt{T-t}}$ with only $s_0$ depending on initial data.  To verify that the blow-up rate agrees with numerical solution we study the quantity $\sqrt{T-t}R(t)$, which should be a linear function of $-\log(T-t)$ with slope independent on initial data.  On the other hand, the shift $-C\,s_0$, should vary with initial data.  In the picture we present results for initial data $u(r,0)=r$ and $u(r,0)=r\pm\sin(r)$, each with its own blow-up time $T$.  The particular values of blow-up time, slope and shift are shown in table \ref{tab:d7-params}.  The blow-up rate $R(t)$ is given by $\partial_r u(0,t)$.}
    \end{figure}

    A more challenging numerical test is to verify the blow-up rate in
    dimension $d=7$.  We expect (cf. equation
    \eqref{eq:generic-k1-d7}) the blow-up rate
    \begin{equation*}
    R(t)=\frac{C\sqrt{T-t}}{(-\log(T-t)-s_0)},\qquad C=\left(\frac{h\gamma}{c_1\, D_1}\right)\frac{1}{\sup_{\xi\ge0} \lvert\frac{dU^*}{d\xi}(\xi)\rvert}.
    \end{equation*}
    This scenario is significantly more difficult to verify than
    \eqref{eq:d8-num-ratio} because in order to see the logarithmic
    correction we must get much closer to the blow-up time $T$.  At
    the same time, the choice of initial data should only influence
    the constant $s_0$, but not $C$.  We start with the relation
    $\partial_r u(0,t)=1/R(t)$, by which we get
    \begin{equation}
    \label{eq:d7-numtest}
    \sqrt{T-t}\,\partial_r u(0,t)=C\,(-\log(T-t)-s_0).
    \end{equation}
    To test our conjectured blow-up rate we plot the left hand side of
    \eqref{eq:d7-numtest} against $-\log(T-t)$, expecting to see a
    linear function after sufficiently long time.  The experimental
    values of $C$, $T$ and $s_0$ are displayed in Table
    \ref{tab:d7-params}, while the relation \eqref{eq:d7-numtest} is
    depicted in Figure \ref{fig:d7}.
\subsection*{Acknowledgments}
\label{sec-2-8}

   We thank Piotr Bizoń for the supervision of this paper.  Special
   thanks are due to Juan L.L. Vel\'azquez and Yukihiro Seki for very
   helpful discussions and suggestions.  This work was supported by a
   Foundation for Polish Science IPP Programme ``Geometry and Topology
   in Physical Models'' and by the NCN Grant No. NN202 030740.

\begin{table}[htb]
\caption{For $d=7$ it holds $\sqrt{T-t}\,\partial_r u(0,t)=C\,(-\log(T-t)-s_0)$, asymptotically as $t\nearrow T$, with $C$ independent of initial data.  In this table we compare values of $C$ and $s_0$ obtained from fitting the asymptotic relation to the numerical solution for various initial data.  The values of $C$ indeed don't change significantly among the tested initial data.} \label{tab:d7-params}
\begin{center}
\begin{tabular}{llll}
 Initial data  &  T           &  $C$        &  $s_0$       \\
\hline
 $r$           &  $0.22913$   &  $0.22512$  &  $-0.43646$  \\
 $r+\sin(r)$   &  $0.066835$  &  $0.22475$  &  $0.45864$   \\
 $r-\sin(r)$   &  $0.44672$   &  $0.22500$  &  $-0.11921$  \\
\hline
\end{tabular}
\end{center}
\end{table}
\section*{Appendix}
\label{sec-3}
\subsection*{Existence and asymptotic form of harmonic maps}
\label{sec-3-1}

\label{sec-3-1-1}

   \begin{theorem}
   \label{th:v-asymp}

   For $d>2+k(2+2\sqrt{2})$, a solution $v(x)$ to equation
   \begin{equation}
   \label{eq:v1}
   v''(x)+(d-2)v'(x)+k(d+k-2)\sin(v)=0
   \end{equation}
   subjected to boundary conditions
   \begin{equation*}
   v(x)=-\pi+2e^{-kx}+\mathcal O(e^{-3kx}),\quad\text{for }x\to-\infty.
   \end{equation*}
   exists and has an asymptotic
   \begin{equation*}
   v(x)=h_+\,e^{-\gamma x}(1+\mathcal O(e^{-2x})+\mathcal O(e^{-\omega x})),\quad\text{for }x\to+\infty
   \end{equation*}
   where $h_+$ is a strictly negative constant, while $\gamma$ and $\omega$ are defined in \eqref{eq:gamma}.
   \end{theorem}

   \emph{Proof.} The proof bases on the analysis of a phase portrait
   spanned by $(v,v')$ of autonomous equation \eqref{eq:v1} and
   consists of three steps.

\begin{description}
\item[Construction of no-escape region] Let us start by defining the
        vector field
        \begin{equation*}
        F(v,v')=(v',-(d-2)v'-k(d+k-2)\sin(v)).
        \end{equation*}
        We are interested in a heteroclinic orbit connecting two
        critical points of $F$, starting at $(-\pi,0)$ and ending at
        $(0,0)$.  We construct a trapping region $\mathcal
        S=\{(v,v')\,|\,-k\sin(v)\le v'\le-\gamma\sin(v),-\pi<v<0\}$,
        which includes critical points $(-\pi,0)$ and $(0,0)$.  No
        integral curve of $F$ starting in $\mathcal S$ can leave
        $\mathcal S$ (see Figure \ref{fig:phase}).

        Indeed, if we define $\underline n(v)=(-k\cos(v),1)$ as a
        normal vector to a curve $v'=-k\sin(v)$, pointing inward of
        $\mathcal S$, by a direct computation we get
        \begin{align*}
        F(v,-k\sin(v))\cdot \underline n(v)&=-k^2\sin(v)(1-\cos(v))
        \end{align*}
        which is positive for $-\pi<v<0$.  Similarly, taking a normal
        vector $\overline n(v)=(-\gamma\cos(v),-1)$ (again directed
        inward $\mathcal S$) to a curve $v'=-\gamma\sin(v)$ gives
        \begin{align*}
        F(v,-\gamma\sin(v))\cdot \overline n(v)&=-\gamma^2\sin(v)(1-\cos(v))
        \end{align*}
        which is also positive for $-\pi<v<0$.  Therefore, the vector
        field $F$ points inward on the whole boundary of $\mathcal S$
        (excluding the stationary points $(0,0)$ and
        $(-\pi,0)$).  This implies that any integral curve of $F$
        starting inside $\mathcal S$ must stay in $\mathcal S$.
\item[Asymptotic of solutions starting in $\mathcal S$] There are two
        stationary points in $\mathcal S$ where a solution can end
        up.  The first one, $(-\pi,0)$, can be ruled out because
        inside $\mathcal S$ vector field $F$ has a nonzero horizontal
        component pointing to the right.  The remaining stationary
        point, $(0,0)$, gives a general asymptotic for of $v$ as
        \begin{equation}
        \label{eq:v1-asymp}
        v(x)=2h_+\cdot e^{x\mu_+}(1+\mathcal O(e^{-2x}))+2h_-\cdot e^{x\mu_-}(1+\mathcal O(e^{-2x}))
        \end{equation}
        where $\mu_\pm<0$ are eigenvalues of $\nabla F(0,0)$
        \begin{equation}
        \label{eq:v1-eigen}
        \mu_+=-\gamma,\quad \mu_-=-\gamma-\omega.
        \end{equation}
        At this point, $h_-$ and $h_+$ are constants depending on
        initial data and there are no restrictions on their
        values.  Because $(v,v')\in\mathcal S$, we have
        $v'<-\gamma\sin(v)<-\gamma v$.  If we combine the latter
        inequality with the asymptotic form of $v$, we get
        \begin{equation}
        \label{eq:v1-inequal}
        -\omega h_-\cdot(1+\mathcal O(e^{-2x}))<0.
        \end{equation}
        On the other hand, from $v<0$ we know that $h_-\cdot(1+\mathcal
        O(e^{-2x}))<0$.  This contradicts with $\omega>0$, so
        $h_+\ne0$.  We can again use $v<0$, this time with leading
        order term proportional to $h_+$, to get $h_+<0$.
\item[Boundary conditions in the thesis guarantee $(v,v')\in\mathcal S$] When
        $x\to-\infty$ the solution $v$ with initial conditions
        $v(x)=-\pi+2e^{kx}+\mathcal O(e^{3kx})$ can be expanded as a
        Taylor series in $e^{x}$ in a following way
        \begin{equation*}
        v(x)=-\pi+2e^{kx}-\frac{2(d+k-2)}{3(d+4k-2)}\,e^{3kx}+\mathcal O(e^{5kx})
        \end{equation*}
        It is a matter of routine computation to show
        that for sufficiently small $x$ we have
        \begin{align}
        \label{eq:v-bound}
        -k\sin(v(x))<v'(x)<-\gamma\sin(v(x)).
        \end{align}
        So $(v,v')\in\mathcal S$ and $v$ has an asymptotic form of
        \eqref{eq:v1-asymp} with $h_+<0$. $\square$
\end{description}
\section*{References}
\label{sec-4}

\bibliographystyle{abbrvnat}
\bibliography{/home/pawel/Documents/Mendeley/library}

\begin{thebibliography}{20}
\providecommand{\natexlab}[1]{#1}
\providecommand{\url}[1]{\texttt{#1}}
\expandafter\ifx\csname urlstyle\endcsname\relax
  \providecommand{\doi}[1]{doi: #1}\else
  \providecommand{\doi}{doi: \begingroup \urlstyle{rm}\Url}\fi

\bibitem[Angenent et~al.(2009)Angenent, Hulshof, and Matano]{Angenent2009}
S.~B. Angenent, J.~Hulshof, and H.~Matano.
\newblock {The Radius of Vanishing Bubbles in Equivariant Harmonic Map Flow
  from \$D\^{}2\$ to \$S\^{}2\$}.
\newblock \emph{SIAM Journal on Mathematical Analysis}, 41\penalty0
  (3):\penalty0 1121--1137, Jan. 2009.
\newblock ISSN 0036-1410.
\newblock \doi{10.1137/070706732}.

\bibitem[Biernat()]{Biernat-MOVCOL}
P.~Biernat.
\newblock {MOVCOL variation in Fortran 95}.
\newblock URL \url{https://github.com/pwl/movcol}.

\bibitem[Biernat and Bizoń(2011)]{Biernat2011}
P.~Biernat and P.~Bizoń.
\newblock {Shrinkers, expanders, and the unique continuation beyond generic
  blowup in the heat flow for harmonic maps between spheres}.
\newblock \emph{Nonlinearity}, 24\penalty0 (8):\penalty0 2211--2228, Aug. 2011.
\newblock ISSN 0951-7715.
\newblock \doi{10.1088/0951-7715/24/8/005}.

\bibitem[Biernat and Seki()]{biernat-seki}
P.~Biernat and Y.~Seki.
\newblock {Type II blow-up mechanisms in harmonic map heat flow, in
  preparation}.

\bibitem[Budd et~al.(2005)Budd, Carretero-Gonz\'{a}lez, and Russell]{Budd2005}
C.~J. Budd, R.~Carretero-Gonz\'{a}lez, and R.~D. Russell.
\newblock {Precise computations of chemotactic collapse using moving mesh
  methods}.
\newblock \emph{Journal of Computational Physics}, 202\penalty0 (2):\penalty0
  463--487, Jan. 2005.
\newblock ISSN 00219991.
\newblock \doi{10.1016/j.jcp.2004.07.010}.

\bibitem[Chen(1989)]{Chen1989}
Y.~Chen.
\newblock {The weak solutions to the evolution problems of harmonic maps}.
\newblock \emph{Mathematische Zeitschrift}, 201\penalty0 (1):\penalty0 69--74,
  Mar. 1989.
\newblock ISSN 0025-5874.
\newblock \doi{10.1007/BF01161995}.

\bibitem[Eells and Sampson(1964)]{Eells1964}
J.~Eells and J.~H. Sampson.
\newblock {Harmonic Mappings of Riemannian Manifolds}.
\newblock \emph{American Journal of Mathematics}, 86\penalty0 (1):\penalty0
  109, Jan. 1964.
\newblock ISSN 00029327.
\newblock \doi{10.2307/2373037}.

\bibitem[Fan(1999)]{Fan1999}
H.~Fan.
\newblock {Existence of the self-similar solutions in the heat flow of harmonic
  maps}.
\newblock \emph{Science in China Series A: Mathematics}, 42\penalty0
  (2):\penalty0 113--132, Feb. 1999.
\newblock ISSN 1006-9283.
\newblock \doi{10.1007/BF02876563}.

\bibitem[Gastel(2002)]{Gastel2002a}
A.~Gastel.
\newblock {Singularities of first kind in the harmonic map and Yang-Mills heat
  flows}.
\newblock \emph{Mathematische Zeitschrift}, 242\penalty0 (1):\penalty0 47--62,
  Feb. 2002.
\newblock ISSN 0025-5874.
\newblock \doi{10.1007/s002090100306}.

\bibitem[Germain and Rupflin(2011)]{Germain2010}
P.~Germain and M.~Rupflin.
\newblock {Selfsimilar expanders of the harmonic map flow}.
\newblock \emph{Annales de l'Institut Henri Poincare (C) Non Linear Analysis},
  28\penalty0 (5):\penalty0 743--773, Sept. 2011.
\newblock ISSN 02941449.
\newblock \doi{10.1016/j.anihpc.2011.06.004}.

\bibitem[Herrero and Vel\'{a}zquez()]{herrero-velazquez-unpublished}
M.~Herrero and J.~J.~L. Vel\'{a}zquez.
\newblock {A blow-up result for semilinear heat equations in the supercritical
  case}.

\bibitem[Herrero and Vel\'{a}zquez(1994)]{Herrero1994}
M.~Herrero and J.~J.~L. Vel\'{a}zquez.
\newblock {Blowup of solutions of supercritical semilinear parabolic
  equations}.
\newblock \emph{C. R. Acad. Sci. Paris S\'{e}r. I Math.}, 2\penalty0
  (3):\penalty0 141--145, 1994.
\newblock URL \url{http://www.sciencedirect.com/science/journal/07644442}.

\bibitem[Herrero and Vel\'{a}zquez(1996)]{Herrero1996}
M.~A. Herrero and J.~J.~L. Vel\'{a}zquez.
\newblock {Singularity patterns in a chemotaxis model}.
\newblock \emph{Mathematische Annalen}, 306\penalty0 (1):\penalty0 583--623,
  Sept. 1996.
\newblock ISSN 0025-5831.
\newblock \doi{10.1007/BF01445268}.

\bibitem[Herrero and Vel\'{a}zquez(1997)]{Herrero1997}
M.~A. Herrero and J.~J.~L. Vel\'{a}zquez.
\newblock {On the Melting of Ice Balls}.
\newblock \emph{SIAM Journal on Mathematical Analysis}, 28\penalty0
  (1):\penalty0 1--32, Jan. 1997.
\newblock ISSN 0036-1410.
\newblock \doi{10.1137/S0036141095282152}.

\bibitem[Huang and Russell()]{Huang-movcol}
W.~Huang and R.~D. Russell.
\newblock {MOVCOL webpage}.
\newblock URL \url{http://www.math.ku.edu/~huang/research/movcol/movcol.html}.

\bibitem[Rapha\"{e}l and Schweyer(2013)]{Raphael2011}
P.~Rapha\"{e}l and R.~Schweyer.
\newblock {Stable Blowup Dynamics for the 1-Corotational Energy Critical
  Harmonic Heat Flow}.
\newblock \emph{Communications on Pure and Applied Mathematics}, 66\penalty0
  (3):\penalty0 414--480, Mar. 2013.
\newblock ISSN 00103640.
\newblock \doi{10.1002/cpa.21435}.

\bibitem[Russell et~al.(2007)Russell, Williams, and Xu]{Russell2007}
R.~D. Russell, J.~F. Williams, and X.~Xu.
\newblock {MOVCOL4: A Moving Mesh Code for Fourth-Order Time-Dependent Partial
  Differential Equations}.
\newblock \emph{SIAM Journal on Scientific Computing}, 29\penalty0
  (1):\penalty0 197--220, Jan. 2007.
\newblock ISSN 1064-8275.
\newblock \doi{10.1137/050643167}.

\bibitem[Struwe(1996)]{Struwe1996}
M.~Struwe.
\newblock {Geometric evolution problems}.
\newblock In \emph{Nonlinear Partial Differential Equations in Differential
  Geometry}, pages 259----339. IAS/Park City Mathematics Series, 1996.

\bibitem[van~den Berg et~al.(2003)van~den Berg, King, and
  Hulshof]{VandenBerg2003}
J.~B. van~den Berg, J.~R. King, and J.~Hulshof.
\newblock {Formal Asymptotics of Bubbling in the Harmonic Map Heat Flow}.
\newblock \emph{SIAM Journal on Applied Mathematics}, 63\penalty0 (5):\penalty0
  1682--1717, Jan. 2003.
\newblock ISSN 0036-1399.
\newblock \doi{10.1137/S0036139902408874}.

\bibitem[Vel\'{a}zquez(2002)]{Vela2002}
J.~J.~L. Vel\'{a}zquez.
\newblock {Stability of Some Mechanisms of Chemotactic Aggregation}.
\newblock \emph{SIAM Journal on Applied Mathematics}, 62\penalty0 (5):\penalty0
  1581, 2002.
\newblock \doi{10.1137/S0036139900380049}.

\end{thebibliography}

\end{document}